\documentclass[11pt,twoside]{amsart}
\date{17 March 2016}
\usepackage{latexsym,amsthm,amsmath,amsfonts,amscd,amssymb}
\usepackage{hyperref}
\usepackage[utf8]{inputenc}
\usepackage[all]{xy}
\usepackage{graphics}
\usepackage{lscape}
\usepackage{array}
\theoremstyle{plain}  
\newtheorem{theorem}{Theorem}[section]

\newtheorem*{theorem*}{Theorem}

\newtheorem{corollary}[theorem]{Corollary}

\newtheorem{proposition}[theorem]{Proposition}

\newtheorem{tech-lemma}[theorem]{Technical Lemma}
\newtheorem{definition}[theorem]{Definition}
\theoremstyle{remark}

\newtheorem{remark}[theorem]{Remark}

\newtheorem*{remark*}{Remark}

\newtheorem*{claim*}{Claim}
\numberwithin{equation}{section}

\renewcommand{\geq}{\geqslant}

\newcommand{\RR}{\mathbb{R}}
\newcommand{\Z}{\mathbb{Z}}
\newcommand{\CC}{\mathbb{C}}

\newcommand{\M}{{\mathcal M}}
\newcommand{\N}{{\mathcal N}}
\newcommand{\cR}{{\mathcal R}}

\newcommand{\vol}{\mathrm{vol}}

\newcommand{\U}{\mathrm{U}}
\newcommand{\GL}{\mathrm{GL}}

\newcommand{\SL}{\mathrm{SL}}

\newcommand{\Sp}{\mathrm{Sp}}

\newcommand{\la}{\langle}
\newcommand{\ra}{\rangle}
\newcommand{\st}{\;|\;}

\DeclareMathOperator{\ad}{ad}
\DeclareMathOperator{\Ad}{Ad}
\DeclareMathOperator{\Aut}{Aut}
\DeclareMathOperator{\Int}{Int}
\DeclareMathOperator{\aut}{aut}

\DeclareMathOperator{\Hom}{Hom}
\DeclareMathOperator{\End}{End}

\newcommand{\liez}{\mathfrak{z}}
\newcommand{\liec}{\mathfrak{c}}
\newcommand{\liep}{\mathfrak{p}}
\newcommand{\liel}{\mathfrak{l}}

\newcommand{\lieh}{\mathfrak{h}}

\newcommand{\lieg}{\mathfrak{g}}
\newcommand{\liegc}{\mathfrak{g}^{\mathbb{C}}}

\let\oldmarginpar\marginpar
\renewcommand\marginpar[1]{\oldmarginpar{\tiny\bf\begin{flushleft} #1
\end{flushleft}}}

\begin{document}

%
%

\title[Connectedness of Higgs bundle moduli for complex reductive Lie groups]
{Connectedness of Higgs bundle moduli for complex reductive Lie groups}
%
%
\author[O. Garc{\'\i}a-Prada]{Oscar Garc{\'\i}a-Prada}
%
\author[A. Oliveira]{André Oliveira}
%
%

\thanks{%
Members of GEAR (Geometric Structures and Representation Varieties).
First  author partially supported by the Ministerio de Econom\'ia y
Competitividad of Spain through Project MTM2010-17717 and Severo Ochoa
Excellence Grant. 
Second author partially supported by CMUP (UID/MAT/00144/2013) and CM-UTAD (PEst-OE/MAT/UI4080/2014), by the Projects PTDC/MAT/120411/2010, PTDC/MAT-GEO/0675/2012 and EXCL/MAT-GEO/0222/2012 and finally by the Post-Doctoral fellowship SFRH/BPD/100996/2014. All these are funded by FCT (Portugal) with national (MEC) and European structural funds (FEDER), under the partnership agreement PT2020.
The second author also thanks the Instituto de Ciencias Matemáticas (ICMAT) in Madrid and the Institute for Mathematical Sciences (IMS) in Singapore --- that visited while preparing the paper --- for the excellent conditions provided.\\
The authors thank Peter Dalakov for useful comments to a preliminary version of the paper and also Tony Pantev for useful discussions concerning the paper \cite{donagi-pantev:2012}.}
\begin{abstract}
We carry an intrinsic approach to the study of the connectedness of the moduli
space $\M_G$ of $G$-Higgs bundles, over a compact Riemann surface, 
when $G$ is a  complex reductive (not necessarily connected) Lie group. 
We prove that the number of connected components of $\M_G$ is indexed by the
corresponding  topological invariants. 
In particular, this gives an alternative proof of the counting by J. Li in
\cite{li:1993} of the number of connected components of the moduli space of 
flat $G$-connections in the case in which  $G$ is connected and semisimple.
\end{abstract}

\maketitle

%
%

 \section{Introduction}\label{section:Introduction}

The topology of the moduli spaces $\M_G$ of $G$-Higgs bundles over a compact Riemann surface $X$ (of genus $g\geq 2$) has been object of intense study in the past decade, mostly by making use of the Morse theoretic techniques introduced by Nigel Hitchin in the seminal paper \cite{hitchin:1987} on Higgs bundles. This procedure uses the fact that the moduli spaces $\M_G$ carry a proper, bounded below real function $f$, from which we can obtain information at least about the connected components, through the study of the subvarieties of local minima of $f$. In some good cases, namely when the spaces $\M_G$ are smooth, the Poincaré polynomial may be calculated, by a study of all critical subvarieties of $f$, since in these cases $f$ is indeed a perfect Morse-Bott function. 

The connected components of $\M_G$ have been object of investigation for many families of real reductive Lie groups $G$, especially after the work of Hitchin \cite{hitchin:1992} where the case $G=\SL(n,\RR)$ was addressed; some references where this subject was studied are \cite{aparicio-garcia-prada:2011,bradlow-garcia-prada-gothen:2003,bradlow-garcia-prada-gothen:2004,bradlow-garcia-prada-gothen:2005,bradlow-garcia-prada-gothen:2013,garcia-prada-gothen-mundet:2013,garcia-prada-mundet:2004,garcia-prada-oliveira:2011,garcia-prada-oliveira:2013,gothen:2001,gothen-oliveira:2012,oliveira:2011}.
The approach has been through a case-by-case study concerning the classes of
$G$. So, the aim of this paper is to take a first step towards the computation
of the number of connected components of $\M_G$, from an intrinsic point of view, in the sense that we do not specify the group $G$. It is an abstract approach using the above mentioned techniques introduced by Hitchin in \cite{hitchin:1987} for $G=\SL(2,\CC)$. Let $c$ represent a topological class of $G$-Higgs bundles and let $\M_G(c)$ denote the subspace of $\M_G$ whose points represent those $G$-Higgs bundles within the class $c$. Of course, $\M_G(c)$ is a union of connected components. We consider a general complex reductive Lie group and we prove the following (see Theorem \ref{thm:main}). 
\begin{theorem}\label{thm:main1}
For any class $c$, the moduli space $\M_G(c)$ is non-empty and connected for any complex reductive Lie group $G$.
\end{theorem}

Some particular cases of this theorem have already been proved in \cite{li:1993,donagi-pantev:2012,lawton-ramras:2014} --- see Remark \ref{rmk:knowncases} --- but all of them use different methods from the ones we use in this paper.

Recall that $\pi_1X$ is a finitely generated group, with $2g$ generators $a_1,b_1,\ldots,a_g,b_g$ such that the product of all commutators $[a_i,b_i]$ is trivial. Now, let $\Gamma$ be the universal central extension of $\pi_1X$, defined as the finitely generated group, with $2g+1$ generators $a_1,b_1,\ldots,a_g,b_g,\delta$, such that $\delta$ lies in the center of $\Gamma$ and $\prod_{i=1}^g[a_i,b_i]=\delta$. Define $\Gamma_\RR$ as $\RR\times_\Z\Gamma$, where $\Z$ is identified as the subgroup of $\Gamma$ generated by $\delta$.
Given any real reductive Lie group $G$, a reductive representation of $\Gamma_\RR$ in $G$ is a continuous homomorphism $\rho:\Gamma_\RR\to G$ which becomes completely reducible when composed with the adjoint representation of $G$. A reductive representation $\rho$ is said to be \emph{central} if $\rho(\RR)$ lies in the centre of $G_0$, the identity component of $G$. Let $\Hom_{\text{cent}}^{\text{red}}(\Gamma_\RR,G)$ be the space of such reductive, central representations. The group $G$ acts on this space by conjugation and we denote by $\cR_G$ the quotient space, usually called the \emph{$G$-character variety of $X$}.

For any real reductive Lie group $G$, non-abelian Hodge theory provides a
homeomorphism between $\M_G$ and $\cR_G$, so our result shows that $\cR_G(c)$
is connected for any complex reductive Lie group and for any class $c$. Here
$\cR_G(c)$ is the subspace of $\cR_G$ whose corresponding central curvature
principal $G$-bundle lies in the topological class $c$ (the homeomorphism
mentioned  above respects the topological classes).
If $G$ is complex, connected and semisimple, then $c$ is
trivial and one consider representations of $\pi_1X$ in $G$. 
In this case, the connectedness of $\cR_G$  is known for more than twenty
years, by the work \cite{li:1993} of Jun Li. Hence, the corresponding result
on the side of Higgs bundles also follows. 
However, Li's methods do not use Higgs bundles. A Higgs bundle approach for
$G=\SL(2,\CC)$ was given by Hitchin in \cite{hitchin:1987} and for
$G=\GL(n,\CC)$ by Simpson \cite{simpson:1994}. In this paper we give  an alternative proof to
the result of Jun Li using Higgs bundles. Moreover, our result is more general
in  the sense that it is valid for reductive and even non-connected complex
Lie groups. Along the way, we prove other results about bundles which we have 
not been able to find  
in the literature. We highlight the following theorem which describes the
stable  and non-simple Higgs bundles, which give rise to orbifold type
singularities  of the moduli space $\M_G$. 
As far as we know, this result does not appear in the literature 
even for principal bundles, and our proof also holds in that case 
(see Theorem \ref{thm:stable and not simple principal bundles}).

\begin{theorem}
Let $G$ be a complex Lie group. For any stable and non-simple $G$-Higgs
bundle, there is a complex reductive Lie subgroup $G'\subset G$ such that the
$G$-Higgs bundle admits a reduction of structure group to $G'$ and it is
stable and simple as a $G'$-Higgs bundle.
\end{theorem}

A natural generalization of this work is to consider a general real reductive
Lie group and we intend to pursue this direction in a different paper. In this
case,  it is well-known that Theorem \ref{thm:main1} does not hold, as there
are some classes of real groups for which $\M_G$ has ``extra'' components.

 \section{$G$-Higgs bundles and topological invariants}\label{section:Higgs bundles}

In this section we introduce the main objects which we shall work with. These are called \emph{$G$-Higgs bundles} and roughly are pairs consisting of a holomorphic bundle and a section of an associated bundle.
$G$-Higgs bundles can be defined on any compact Kähler manifold (cf. \cite{simpson:1992}), and $G$ may be any real reductive Lie group (see for example \cite{bradlow-garcia-prada-gothen:2005}), but we will restrict ourselves to $G$-Higgs bundles over compact Riemann surfaces, and such that $G$ is a complex reductive Lie group.

Fix a compact and connected Riemann surface $X$ of genus $g\geq 2$. Let $K=T^*X^{1,0}$ be its canonical line bundle. Let $G$ be a complex reductive Lie group.
Given a principal $G$-bundle $E_G$, denote by $\ad(E_G)$ the adjoint bundle of $E_G$, that is the vector bundle obtained from $E_G$ under the adjoint representation of $G$ on its Lie algebra $\lieg$: $$\ad(E_G)=E_G\times_G\lieg.$$

\begin{definition}\label{definition of Higgs bundle}
A \emph{$G$-Higgs bundle} over $X$ is a pair
$(E_G,\varphi)$ where $E_G$ is a holomorphic principal $G$-bundle over
$X$ and $\varphi$ is a holomorphic section of $\ad(E_G)\otimes K$. The section $\varphi$ is usually called the \emph{Higgs field}.
\end{definition}

As an example, a $\GL(n,\CC)$-Higgs bundle or simply a Higgs bundle is, in terms of holomorphic vector bundles, a pair $(V,\varphi)$ with $V$ a holomorphic rank $n$ vector bundle and $\varphi\in H^0(X,\End(V)\otimes K)$, whereas for $\SL(n,\CC)$-Higgs bundles, $V$ is required to have trivial determinant and $\varphi$ must be trace-free.
Higgs bundles were first introduced by Hitchin in \cite{hitchin:1987}, for
$G=\SL(2,\CC)$, while studying the 
self-duality equations (now known as Hitchin equations) on Riemann surfaces.

The topological class of a $G$-Higgs bundle is given by the topological class of the underlying $G$-principal bundle. If $G$ is connected, the topological classification of principal $G$-bundles over the compact Riemann surface $X$ is well-known to be given by elements of $\pi_1G$ (cf. \cite[Proposition 5.1]{ramanathan:1975}). For a not necessarily connected group $G$, the topological classification of G-bundles is more subtle, and we only briefly sketch it; details can be found in \cite[\S 2]{oliveira:2008} and in \cite[Prop. 3.1]{oliveira:2011}. 
First we assume that $\pi_0G$ is an abelian group; this assumption is only needed for Theorem \ref{thm:topclass} below and nothing else.
Given a principal $G$-bundle $E_G$ on $X$, let $m_1(E_G)$ be the induced (flat) $\pi_0G$-bundle. 
This gives a first topological invariant of $E_G$, as $$m_1(E_G)\in H^1(X,\pi_0G)\cong(\pi_0G)^{2g}.$$
Now, $\pi_0G$ acts on $\pi_1G$ through the conjugation action of $G$ on itself. Fix a class $m_1\in H^1(X,\pi_0G)$ and let $\pi_1\mathcal{G}_{m_1}$ be the flat $\pi_1G$-bundle associated to $m_1: \pi_1X\to\pi_0G$ via the action $\pi_0G\to\Aut(\pi_1G)$.
We can consider cohomology with values in the local coefficient system $\pi_1\mathcal{G}_{m_1}$. 
In fact, $\pi_0G$ also acts on $H^2(X,\pi_1\mathcal{G}_{m_1})$ through $\pi_0G\to\Aut(\pi_1G)$, and the next result says that topological $G$-bundles on $X$ with the first class $m_1$ fixed are classified by elements in the quotient space.

\begin{theorem}[\cite{oliveira:2011}, Proposition 3.1]\label{thm:topclass}
Let $G$ be a Lie group with $\pi_0G$ abelian.
There is a bijection between the set of isomorphism classes of continuous principal $G$-bundles over the surface $X$ with invariant $m_1\in (\pi_0G)^{2g}$ and the quotient set $H^2(X,\pi_1\mathcal{G}_{m_1})/\pi_0G$.
 \end{theorem}

In fact, this theorem is valid not only on surfaces, but on any $2$-dimensional connected $CW$-complex. 
Observe that if $G$ is connected, the preceding theorem gives the bijection between topological classes of $G$-Higgs bundles over $X$ and $H^2(X,\pi_1G)\cong\pi_1G$, as we already mentioned.

\section{Semistability and moduli spaces}

\subsection{Semistability}

Let $G$ be a complex reductive Lie group. In order to consider moduli spaces of $G$-Higgs bundles we need the corresponding notions of (semi,poly)stability.
We briefly recall the main definitions. The main reference is \cite{garcia-prada-gothen-mundet:2008}, where all these general notions are deduced in detail.

We recall first some definitions. Let $\lieg$ be the Lie algebra of $G$ and $\liez$ its centre. Then $\lieg=\liez\oplus\lieg_{ss}$, where $\lieg_{ss}=[\lieg,\lieg]$ is the semisimple part of $\lieg$. Given a Cartan subalgebra $\liec$ of $\lieg_{ss}$, we will consider roots of $\lieg$ as forms on $\liec$ extended by zero on $\liez$.
Let $R$ be the set of such roots and for $\alpha\in R$, let $\lieg_\alpha$ be the corresponding root space, so that we have the corresponding decomposition: 
$$\lieg=\liez\oplus\liec\oplus\bigoplus_{\alpha\in R}\lieg_\alpha.$$ Let $\Delta\subset R$ be the system of simple roots.


Let $\lieh$ be the Lie algebra of a maximal compact subgroup of $G$. Given $s\in i\lieh$,
\begin{equation}\label{eq:parabolicPs}
P_s=\{g\in G\st e^{ts}ge^{-ts}\text{ remains bounded  as }t\to\infty\},
\end{equation}
is a parabolic subgroup of $G$, whose corresponding parabolic subalgebra of $\lieg$ is $$\liep_s=\{v\in\lieg\st\Ad(e^{ts})(v)\text{ remains bounded  as }t\to\infty\}.$$
If, moreover, we define 
\begin{equation}\label{eq:LeviLs}
L_s=\{g\in G\st\lim_{t\to\infty} e^{ts}ge^{-ts}=g\}
\end{equation}
then $L_s\subset P_s$ is a Levi subgroup of $P_s$, and 
$$\liel_s=\{v\in\lieg\st\lim_{t\to\infty}\Ad(e^{ts})(v)=0\}$$ is the corresponding Levi subalgebra of $\liep_s$.

In case $G$ is connected, every parabolic subgroup $P$ is of the form \eqref{eq:parabolicPs} for some $s\in i\lieh$; the same holds for the Levi subgroups. For $G$ non-connected that may not be the case (cf. \cite[Remark 5.3]{martin:2003}). However, in order to define semistability, the parabolic subgroups which need to be considered are precisely the ones of the form \eqref{eq:parabolicPs}. Hence, for simplicity, and when no explicit mention to $s\in i\lieh$ is needed, we refer to these as the parabolic subgroups of $G$, keeping in mind that we mean the groups defined by \eqref{eq:parabolicPs}. We will do the same for the Levi subgroups, referring to \eqref{eq:LeviLs}.

Let then $P$ be a parabolic subgroup of $G$. A character of the Lie algebra
$\liep$ of $P$ is a complex linear map $\liep\to\CC$ which factors through
$\liep/[\liep,\liep]$. Let $\liel\subset\liep$ be the corresponding Levi subalgebra and let $\liez_\liel$ be the centre of $\liel$. Then, one has that $(\liep/[\liep,\liep])^*\cong\liez_\liel^*$, so the characters of $\liep$ are indeed classified by elements of $\liez_\liel^*$.
Since $\lieg$ is reductive, the Killing form on its semisimple part extends to a non-degenerate invariant $\CC$-bilinear pairing $\la\cdot,\cdot\ra$ on $\lieg$, which yields an isomorphism $\liez_\liel^*\cong\liez_\liel$.
Thus, a character $\chi_*\in\liez_\liel^*$ of $\liep$ is uniquely determined by an element $s_{\chi_*}\in\liez_\liel$. Indeed, it can be shown that $\liez_\liel\subset i\lieh$, so that $s_{\chi_*}\in i\lieh$.
Now, the character $\chi_*$ of $\liep$ is said to be \emph{antidominant} if $\liep\subset\liep_{s_{\chi_*}}$ and  \emph{strictly antidominant} if $\liep=\liep_{s_{\chi_*}}$. 

Given a character $\chi:P\to\CC^*$ of $P$, denote by $\chi_*$ the corresponding character of $\liep$. We say that $\chi$ is \emph{(strictly) antidominant} if $\chi_*$ is.

Let $E_G$ be a holomorphic principal $G$-bundle on  $X$ and let $P$ be a
parabolic subgroup of $G$. Denote by $E_G(G/P)$ the holomorphic bundle with
fibre $G/P$ associated to $E_G$ and to the standard action of $G$ on
$G/P$. The bundle $E_G(G/P)$ is canonically isomorphic $E_G/P$. Let $\sigma\in H^0(X,E_G/P)$, that is, a reduction of the structure group of $E_G$ to $P$, and denote by $E_{\sigma}\subset E_G$ the corresponding holomorphic principal $P$-bundle on $X$. So, $E_\sigma$ is the pullback of the principal $P$-bundle $E_G\to E_G/P$ under $\sigma:X\to E_G/P$.
Given the holomorphic principal $P$-bundle $E_\sigma$, we consider the adjoint bundle $$\ad(E_\sigma)=E_\sigma\times_{P}\liep$$ and if we have a further reduction of structure group $\sigma_{L}$ of $E_\sigma$ to a principal $L$-bundle $E_{\sigma_{L}}$, then we can consider also $$\ad(E_{\sigma_{L}})=E_{\sigma_{L}}\times_{L}\liel.$$

Let $\chi:P\to\CC^*$ be an antidominant character of $P$.
The \emph{degree} of
$E_\sigma$, with respect to $\chi$, denoted by $\deg_{\chi}(E_\sigma)$, is the degree of the line bundle obtained by extending the structure group of $E_\sigma$ through $\chi$. In other words,
\begin{equation}\label{def:degree}
\deg_{\chi}(E_\sigma)=\deg(E_\sigma\times_\chi\CC^*).
\end{equation}

Here is the general definition of (semi,poly)stability of $G$-Higgs bundles over a compact Riemann surface $X$ and for $G$ a reductive complex Lie group. It depends on a parameter $\alpha\in i\liez_\lieh$, where $\liez_\lieh$ denotes the center of $\lieh$. Recall that $\la\cdot,\cdot\ra$ denotes an invariant $\CC$-bilinear pairing on
$\lieg$ extending the Killing form on the semisimple part $\lieg_{ss}$. Details may be found in \cite{garcia-prada-gothen-mundet:2008}, where these conditions were defined in the more general setting of any principal pairs for any real reductive Lie group.

\begin{definition}\label{def:semipoly}
Let $\alpha\in i\liez_\lieh$. A $G$-Higgs bundle $(E_G,\varphi)$ over $X$ is:
\begin{itemize}
\item \emph{$\alpha$-semistable} if $\deg_{\chi}(E_\sigma)-\la\alpha,s_{\chi_*}\ra\geq 0$, for any parabolic subgroup $P$ of $G$, any non-trivial antidominant character $\chi$ of $P$ and any reduction of structure group $\sigma$ of $E_G$ to $P$ such that $\varphi\in H^0(X,\ad(E_\sigma)\otimes K)$.
\item \emph{$\alpha$-stable} if $\deg_{\chi}(E_\sigma)-\la\alpha,s_{\chi_*}\ra> 0$, for any parabolic subgroup $P$ of $G$, any non-trivial antidominant character $\chi$ of $P$ and any reduction of structure group $\sigma$ of $E_G$ to $P$ such that $\varphi\in H^0(X,\ad(E_\sigma)\otimes K)$.
\item \emph{$\alpha$-polystable} if it is $\alpha$-semistable and if the following holds. Suppose that $\deg_{\chi}(E_\sigma)-\la\alpha,s_{\chi_*}\ra=0$, for some parabolic subgroup $P\subset G$, some non-trivial strictly antidominant character $\chi$ of $P$ and some reduction of structure group $\sigma$ of $E_G$ to $P$ such that $\varphi\in H^0(X,\ad(E_\sigma)\otimes K)$. Then there is a further holomorphic reduction of structure group $\sigma_L$ of $E_\sigma$ to the Levi subgroup $L$ of $P$ such that $\varphi\in H^0(X,\ad(E_{\sigma_L})\otimes K)$.
\end{itemize}
\end{definition}

\begin{remark}\label{rm:(semi)stability of bundles}
A $G$-Higgs bundle with $\varphi=0$ is a holomorphic principal $G$-bundle and a (semi)stability condition for these objects over compact Riemann surfaces was established by Ramanathan in \cite{ramanathan:1975}. A direct generalization of Ramanathan's condition to the $G$-Higgs bundle case (for $G$ complex) is given in 
\cite{biswas-gomez:2008}. In both cases, we see that the stability does not depend on any parameter $\alpha$. There is however no discrepancy between both (semi)stability conditions, because in \cite{ramanathan:1975,biswas-gomez:2008} the authors only consider characters which are trivial on the center of $G$. The corresponding ones on the Lie algebra are thus orthogonal to $\alpha$ with respect to the pairing $\la\cdot,\cdot\ra$, hence the parameter vanishes on the conditions. The above  definition of \cite{garcia-prada-gothen-mundet:2008} is finer in the sense that it makes precise that there is a parameter involved.
One can say that the precise relation between both conditions is hence that a $G$-Higgs bundle $(E_G,\varphi)$ is (semi)stable in the sense of \cite{ramanathan:1975,biswas-gomez:2008} if and only if it is $\alpha$-(semi)stable in the sense of \cite{garcia-prada-gothen-mundet:2008}, for some $\alpha$.
The significance of the parameter $\alpha$ is more obvious in the
generalization of the notions of Higgs bundle and stability for real reductive
Lie groups (see \cite{garcia-prada-gothen-mundet:2008}). When $G$ is complex,
given a $G$-Higgs bundle $(E_G,\varphi)$, 
we shall see in fact below (cf. Proposition \ref{prop:fixed parameter}) that
the value of the parameter $\alpha$ is uniquely determined by the topological type of $
(E_G,\varphi)$. In other words, $(E_G,\varphi)$ can only be
$\alpha$-polystable  if $\alpha$ is the element in $ i\liez_\lieh$ determined by the topological type of $E_G$.
\end{remark}

Denote by $\M^\alpha_G(c)$ the moduli space of $\alpha$-semistable $G$-Higgs
bundles with fixed topological class $c$ over the Riemann surface $X$. As
usual, the moduli space $\M^\alpha_G(c)$ can also be viewed as parametrizing
isomorphism classes of $\alpha$-polystable $G$-Higgs bundles. The moduli space
$\M^\alpha_G(c)$ has the structure of a quasi-projective variety, as one can
see  from the Schmitt's general Geometric Invariant Theory construction of 
moduli of decorated bundles (cf. \cite{schmitt:2008}), which applies in 
particular to the case of $G$-Higgs bundles, without assuming the
connectedness  of $G$ (cf. \cite[Remark 2.7.5.4]{schmitt:2008}). For related 
constructions also for  higher dimensional projective varieties one can look 
at the work of
Simpson \cite{simpson:1994,simpson:1995}

\subsection{Hitchin equations and $\alpha$-polystability condition}\label{section:Hitchin equations}

Let $(E_G,\varphi)$ be a $G$-Higgs bundle over $X$. By an abuse of notation, we shall denote the $C^\infty$-objects underlying $E_G$ and $\varphi$ by the same symbols. Then the Higgs field may be viewed as a $(1, 0)$-form on $X$ with values in $\ad(E_G)$, $\varphi\in \Omega^{1,0}(X,\ad(E_G))$. Let $H\subset G$ be a maximal compact subgroup. Then its Lie algebra $\lieh$ is a compact form of $\lieg$.
Given a $C^\infty$ reduction of structure group $h$ of $E_G$ to $H$, let $F_h$
be the curvature of the corresponding Chern connection (the unique $H$-connection compatible with the holomorphic structure of $E_G$). Let also
$\tau_h:\Omega^{1,0}(X,\ad(E_G))\to \Omega^{0,1}(X,\ad(E_G))$ be the
involution given by the combination of complex conjugation on complex
$1$-forms with the compact conjugation on $\liegc$ which determines the
compact form $\lieh$, and which is given fibrewise by the metric $h$.  
Let $\omega$ be a volume form of $X$.
The Hitchin-Kobayashi correspondence asserts the following.

\begin{theorem}\label{thm:Hit-Kob}
A $G$-Higgs bundle $(E_G,\varphi)$ is $\alpha$-polystable if and only if there is a reduction of structure group $h$ of $E_G$ from $G$ to $H$ that satisfies the \emph{Hitchin equation} $F_h-[\varphi,\tau_h(\varphi)]=-i \alpha\omega$.
\end{theorem}

A proof of this correspondence can be found in \cite[Theorems 2.24 and 3.21]{garcia-prada-gothen-mundet:2008} (see also \cite{bradlow-garcia-prada-mundet:2003}), in fact in a much more general setting than the one we are considering here. Indeed, this correspondence was first proved for $G=\SL(2,\CC)$ by Hitchin in \cite{hitchin:1987}.

The polystability condition for a $G$-Higgs bundle $(E_G,\varphi)$ depends, in principle, of a parameter $\alpha\in i\liez_\lieh$, but as $G$ is complex, $\alpha$ is indeed fixed by the topological type of $E_G$, as we now show.

\begin{proposition}\label{prop:fixed parameter}
Let $(E_G,\varphi)$ be a $G$-Higgs bundle of topological type $c$. Then there is a unique value of $\alpha\in i\liez_\lieh$, determined by $c$, for which $(E_G,\varphi)$ can be $\alpha$-polystable.
\end{proposition}
\proof
This can be seen combining Chern-Weil theory and Theorem \ref{thm:Hit-Kob}. Indeed, if $(E_G,\varphi)$ is $\alpha$-polystable then it corresponds to a solution $h$ of 
$F_h-[\varphi,\tau_h(\varphi)]=-i\alpha\omega$. By applying any degree one $H$-invariant polynomial $p$ on $\lieh$ to the equation and integrating over $X$, we obtain
\begin{equation}\label{eq:obtaining alpha}
p(\alpha)=\frac{i}{\vol(X)}\int_X[p(F_h)].
\end{equation} Notice that we have used here that $G$ is a complex group because in this case $[\varphi,\tau_h(\varphi)]$ is in the semisimple part of $\lieh$, hence any $p$ vanishes on it. Chern-Weil theory implies that the cohomology class $[p(F_h)]\in H^2(X,\CC)$ only depends on the topological class $c$ of $E_G$. Since a degree one $H$-invariant polynomial of $\lieh$ is a linear map $p:\lieh\to\CC$ which factors through $\lieh/[\lieh,\lieh]\cong\liez_\lieh$, the space of such polynomials is identified with the dual of $\liez_\lieh$. As
$\alpha\in i\liez_\lieh$, then applying \eqref{eq:obtaining alpha} simultaneously for a basis of $\liez_\lieh^*$ determines $\alpha$, as required.
\endproof 


For example, as mentioned above, a $\GL(n,\CC)$-Higgs bundle is equivalent
to  a pair $(V,\varphi)$ where $V$ is a rank $n$ holomorphic vector bundle and $\varphi\in H^0(\End(V)\otimes K)$. The topological type of $(V,\varphi)$ is given by the degree $d$ of $V$. If we normalize the volume of $X$ to be $\vol(X)=2\pi$ and take the trace as a base of the space of degree one invariant polynomials in $\mathfrak{u}(n)$, then $(V,\varphi)$ can only be $\alpha$-(semi,poly)stable in the sense of Definition \ref{def:semipoly} if $\alpha$ equals the slope of $V$, i.e., $\alpha=d/n$. In this case, then one checks that indeed the $d/n$-(semi,poly)stability condition is equivalent to the usual (semi,poly)stability condition comparing the slopes of $V$ and of its $\varphi$-invariants subbundles.

Fixing the topological type $c$ thus fixes $\alpha$, and  hence 
we can  just write $\M_G(c)$ for the moduli space of polystable 
$G$-Higgs bundles over $X$, of topological type $c$, where it is 
implicit that we are using the $\alpha$ given by $c$.


Let \begin{equation}\label{eq:modulispace}
\M_G=\bigsqcup_c\M_G(c),
\end{equation}
where $c$ runs over all possible topological types of $G$-Higgs bundles, according to Theorem \ref{thm:topclass} (if $\pi_0G$ abelian).

 \subsection{Non-emptiness}\label{section:nonemptiness}

We now want to prove that the spaces $\M_G(c)$ are non-empty. For that it is enough to prove the existence, for any $c$, of a polystable $G$-principal bundle with topological type $c$ (so a $G$-Higgs bundle with vanishing Higgs field). This is well-known to be true in case $G$ is connected (see the first part of the proof of Theorem \ref{thm:nonemptiness} below), so the main purpose is to prove non-emptiness of $\M_G(c)$ for non-connected $G$.

Let $H\subset G$ be a maximal compact subgroup. Let $E_H$ be a $C^\infty$ $H$-principal bundle over $X$, with projection map \[p_X:E_H\to X.\] Let $H_0$ be the connected component of the identity of $H$, so that we have 
\begin{equation}\label{eq:sespi0(G)}
1\to H_0\to H\xrightarrow{\pi} \pi_0H\to 1.
\end{equation}
Let $Y$ be the quotient $E_H/H_0$. Then the projection $E_H\to E_H/H_0$ is a $C^\infty$ $H_0$-principal bundle over $Y$. Denote this $H_0$-bundle by 
\begin{equation}\label{eq:H_0-bundletoY}
p_Y:E_{H_0}\to Y.
\end{equation} Clearly the total spaces of $E_H$ and of $E_{H_0}$ are the same; they just project to different basis, having hence different fibres and structure groups. Notice also that $p:Y\to X$ is an unramified covering. It is in fact a $\pi_0H$-principal bundle over $X$, and $p_X=p\circ p_Y$.

\begin{proposition}
Any connection $A_0$ on $E_{H_0}\to Y$ induces naturally a connection $A$ on $E_H\to X$. 
\end{proposition}
\proof
Let $A_0$ be a connection on $E_{H_0}$. Then $A_0\in\Omega^1(E_{H_0},\lieh)$ is an $\lieh$-valued $1$-form on $E_{H_0}$, such that for every $x\in E_{H_0}$, $A_{0,x}:T_xE_{H_0}\to\lieh$ is a splitting of the exact sequence 
\[0\to\lieh\xrightarrow{v_{x,*}^0} T_xE_{H_0}\xrightarrow{(p_{Y,*})_x} T_{p_Y(x)}Y\to 0\] where $v_{x,*}^0$ is the differential at the identity of the map $v_x^0:H_0\to E_{H_0}$, $h\mapsto x\cdot h$, given by the right $H_0$-action on $E_{H_0}$.
Moreover, if $H_x=\ker A_{0,x}\subset T_xE_{H_0}$ is the horizontal subspace of $T_xE_{H_0}$, then $(R_h)_{*,x}(H_x)=H_{x\cdot h}$ for any $h\in H_0$, where $R_h:E_{H_0}\to E_{H_0}$, $x\mapsto x\cdot h$. Recall that $H_x$ is isomorphic to $T_{p_Y(x)}Y$ via $(p_{Y,*})_x$.
 
Since the total spaces $E_{H_0}$ and $E_H$ are the same, we can define $A\in\Omega^1(E_H,\lieh)$ as $A=A_0$. We have to see that $A$ is indeed a connection on $E_H$.

For every $x\in E_H$, it is clear that $A_x:T_xE_H\to\lieh$ is a splitting of
\[0\to\lieh\xrightarrow{v_{x,*}} T_xE_H\xrightarrow{(p_{X,*})_x} T_{p_X(x)}X\to 0,\] where $v_{x,*}$ is the differential at the identity of $v_x:H\to E_H$, $h\mapsto x\cdot h$. Note that this restricts to the given $H_0$-action on $E_{H_0}$, so $v_{x,*}=v_{x,*}^0$.
Note also that 
\begin{equation}\label{eq:isom tangent spaces}
T_{p_X(x)}X\cong T_{p_Y(x)}Y,
\end{equation} for every $x\in E_H$, since $Y\to X$ is an unramified covering.

Pick any $x\in E_H$ and any $h\in H$. Let $H_x=\ker A_x$. We want to see $(R_h)_{*,x}(H_x)=H_{x\cdot h}$. 
Clearly, the right $H$-action on $E_H$ is equivariant with the right $\pi_0H$-action on $Y$, so we have the following commutative diagram, where we are implicitly using \eqref{eq:isom tangent spaces}
\[\xymatrix{0\ar[r]&\lieh\ar[r]^(.4){v_{x\cdot h}}&T_{x\cdot h}E_H\ar[r]^(.45){(p_{Y,*})_{x\cdot h}}&T_{p_Y(x\cdot h)}Y\ar[r]&0\\
0\ar[r]&\lieh\ar[r]^(.4){v_x}\ar[u]^{=}&T_xE_H\ar[r]^(.45){(p_{Y,*})_x}\ar[u]^{(R_h)_{*,x}}&T_{p_Y(x)}Y\ar[r]\ar[u]_{(R_{\pi(g)})_{*,x}}&0}\]
and where $\pi(h)$ is the projection of $h\in H$ in $\pi_0H$ as in \eqref{eq:sespi0(G)}.
Since $(p_{Y,*})_x$ and $(p_{Y,*})_{x\cdot h}$ map $H_x$ and $H_{x\cdot h}$ respectively onto $T_{p_Y(x)}Y$ and $T_{p_Y(x\cdot h)}Y$, respectively, we conclude that indeed $(R_h)_{*,x}(H_x)=H_{x\cdot h}$. 

Hence, $A$ is a connection on $E_H\to X$.
\endproof

\begin{theorem}\label{thm:nonemptiness}
For any topological type $c$ given by Theorem \ref{thm:topclass}, there exists a polystable $G$-principal bundle.
\end{theorem}
\proof
Take $H$ as above. If $G$ is connected (then so is $H$), from Proposition 6.16 of \cite{atiyah-bott:1982} it follows that given a $C^\infty$ $H$-principal bundle $E_H\to X$ with any topological type $c\in\pi_1G$, it admits a Hermitian-Einstein connection (that is a connection $A$ whose curvature is constant and defined by an element in $\liez_\lieh$, the centre of the Lie algebra of $H$). From \cite{ramanathan-subramanian:1988}, that implies the polystability of the holomorphic $G$-principal bundle $E_G$ associated to $E_H$ and with the holomorphic structure such that $A$ is the Chern connection on $E_G$ (so the holomorphic structure on $E_G$ is given by $\bar\partial_A=A^{0,1}$, the $(0,1)$-part of $A$).

Assume that $G$ (hence $H$) is not connected. Let $H_0$ be the component of the identity.
Take an $H$-principal bundle $E_H\to X$ in the $C^\infty$ category, with any topological type $c$ (given by Theorem \ref{thm:topclass}).  From it, construct the $C^\infty$ $H_0$-principal bundle $E_{H_0}\to Y=E_H/H_0$ as in \eqref{eq:H_0-bundletoY}.
Since $H_0$ is connected, $E_{H_0}$ admits a Hermitian-Einstein connection $A_0\in\Omega^1(E_{H_0},\lieh)$ such that the holomorphic $G_0$-principal bundle $E_{G_0}\to Y$ associated to $E_{H_0}$ and to $\bar\partial_{A_0}$ is polystable. 
From the preceding proposition, $A_0$ yields naturally a connection $A$ on $E_H$. Take the corresponding holomorphic $G$-principal bundle $E_G\to X$, coming from $\bar\partial_A$. The topological type of $E_G$ is still the given $c$. Finally, since $E_{G_0}$ is polystable, then from \cite{biswas-gomez:2014}, so is $E_G$.
\endproof

\begin{corollary}\label{cor:Mc(G)non-empty}
For any topological type $c$ given by Theorem \ref{thm:topclass}, $\M_G(c)$ is non-empty.
\end{corollary}

 \subsection{Deformation theory}\label{section:def theory}

We briefly recall  the deformation theory of $G$-Higgs bundles and, in particular, the identification of the tangent space of $\M_G$ at the smooth points. Details can be found for instance in \cite{garcia-prada-gothen-mundet:2008,bradlow-garcia-prada-gothen:2013}.

\begin{definition}
Let $(E_G,\varphi)$ a $G$-Higgs bundle over $X$. The \emph{deformation complex} of $(E_G,\varphi)$ is the complex of sheaves on $X$ given by
\begin{equation}\label{eq:deformation complex}
C^\bullet(E_G,\varphi):\ad(E_G)\xrightarrow{\text{ad}(\varphi)}\ad(E_G)\otimes K.
\end{equation}
\end{definition}

\begin{proposition}\label{prop:deformation for Higgs}
Let $(E_G,\varphi)$ be a $G$-Higgs bundle over $X$.
\begin{enumerate}
    \item The
infinitesimal deformation space of $(E_G,\varphi)$ is isomorphic to the
first hypercohomology group $\mathbb{H}^1(C^\bullet(E_G,\varphi))$ of the
complex $C^\bullet(E_G,\varphi)$;
    \item There is a long exact sequence
\begin{equation*}
\begin{split}
0&\longrightarrow\mathbb{H}^0(C^\bullet(E_G,\varphi))\longrightarrow
H^0(\ad(E_G))\longrightarrow
H^0(\ad(E_G)\otimes
K)\longrightarrow\\
&\longrightarrow\mathbb{H}^1(C^\bullet(E_G,\varphi))\longrightarrow H^1(\ad(E_G))\longrightarrow
H^1(\ad(E_G)\otimes K)\longrightarrow\\
&\longrightarrow\mathbb{H}^2(C^\bullet(E_G,\varphi))\longrightarrow 0
\end{split}
\end{equation*}
where the maps $H^i(\ad(E_G))\to H^i(\ad(E_G)\otimes K)$ are induced by
$\mathrm{ad}(\varphi)$.
\end{enumerate}
\end{proposition}

In particular, it follows from (i) of the proposition that if $(E_G,\varphi)$ represents a smooth point of the moduli space $\M_G(c)$, then $\mathbb{H}^1(C^\bullet(E_G,\varphi))$ is canonically isomorphic to the tangent space at this point.
From this one has that 
\begin{equation}\label{eq:dimH1}
\begin{split}
\dim\mathbb{H}^1(C^\bullet(E_G,\varphi))&=\chi(\ad(E_G)\otimes K)-\chi(\ad(E_G))+\\
&+\dim\mathbb{H}^0(C^\bullet(E_G,\varphi))+\dim\mathbb{H}^2(C^\bullet(E_G,\varphi))
\end{split}
\end{equation}
where $\chi=\dim H^0-\dim H^1$ denotes the Euler characteristic. 

Let $\Aut(E_G,\varphi)$ be the  group of automorphisms of $(E_G,\varphi)$. 
\begin{equation}\label{eq:autom}
\Aut(E_G,\varphi)=\{s\in\Aut(E_G)\st\Ad(s)(\varphi)=\varphi\},
\end{equation} where we recall that $\Aut(E_G)=H^0(\Ad(E_G))$ and $\Ad(E_G)=E_G\times_GG$ with $G$ acting by conjugation. Let also $\aut(E_G,\varphi)$ be the space of infinitesimal automorphisms of $(E_G,\varphi)$, defined as
\begin{equation}\label{eq:infinit autom}
\aut(E_G,\varphi)=\{s\in\aut(E_G))\st\ad(s)(\varphi)=0\},
\end{equation} with $\aut(E_G)=H^0(\ad(E_G))$ Clearly, from (ii) of Proposition \ref{prop:deformation for Higgs}, we have $$\mathbb{H}^0(C^\bullet(E_G,\varphi))\cong\aut(E_G,\varphi).$$

\begin{remark}\label{rem:identif-Aut-G}
Let $f\in\Aut(E_G)$, so that $f$ is a global section of $\Ad(E_G)$. Given $x\in X$, the element $f(x)\in\Ad(E_G)_x$ is identified with a conjugacy class of an element of $g\in G$ (via the identification of $\Ad(E_G)_x$ with $G$, up to an inner automorphism of $G$). Now, the \emph{closure} of the orbit of $g$ under conjugation is independent of the point $x$. This is because the value at $g$ of the $G$-invariant polynomials determine the closure of the orbit of $g$ under conjugation. Since the coefficients of these polynomials are holomorphic functions and $X$ is compact, they are constants, so the closure of the conjugation orbit of $g$ is independent of $x\in X$. 
Indeed, in the cases we will deal with, an element of $f\in\Aut(E_G)$ will determine really a conjugation class $[g]$ of an element $g\in G$. This will follow because the orbits, under conjugation, of such $g\in G$ will be closed.
\end{remark}

Let $Z(G)$ denote the centre of $G$. By Remark \ref{rem:identif-Aut-G}, and since $Z(G)$ acts trivially by conjugation, we can consider $Z(G)$ as a subgroup of $\Aut(E,\varphi)$. If $\liez$ denotes the Lie algebra of $Z(G)$ then, analogously, $\liez$ is a subalgebra of $\aut(E_G,\varphi)$.

\begin{definition}
A $G$-Higgs bundle is \emph{simple} if $\Aut(E_G,\varphi)\cong Z(G)$.
\end{definition}

In order for  $(E_G,\varphi)$ to represent a smooth point  of the moduli space $\M_G$, $\dim\mathbb{H}^0(C^\bullet(E_G,\varphi))$ and $\dim\mathbb{H}^2(C^\bullet(E_G,\varphi))$ must have the minimum possible value. 
Indeed, Serre duality provides an isomorphism
$\mathbb{H}^2(C^\bullet(E_G,\varphi))\cong\mathbb{H}^0(C^\bullet(E_G,\varphi))^*$, and we have the following result (cf. \cite{bradlow-garcia-prada-gothen:2013}):

\begin{proposition}\label{prop:stable and simple imply smooth}
Let $(E_G,\varphi)$ be a stable and simple $G$-Higgs bundle. Then it represents a smooth point of the moduli space $\M_G(c)$.
\end{proposition}

If we are in the situation of the previous proposition, then the dimension of component of the moduli space containing $(E_G,\varphi)$ equals the dimension of $\mathbb{H}^1(C^\bullet(E_G,\varphi))$ which, from \eqref{eq:dimH1}, becomes $\dim(G)(2g-2)+2\dim Z(G)$. Notice that this is twice the dimension of the moduli space $\N_G$ of principal $G$-bundles on $X$. Indeed, by considering the cotangent bundle of $\N_G$, one is naturally lead to $G$-Higgs bundles, for $G$ complex. In fact, $\M_G$ strictly contains the cotangent bundle of $\N_G$.

 \subsection{Stable and not simple $G$-Higgs bundles}\label{section:stable and not simple}

Stable and non-simple $G$-Higgs bundles are represented by points of the
moduli space which may be orbifold type singularities. In this section we show
that such $G$-Higgs bundles always reduce to a stable and simple $G'$-Higgs
bundle for some smaller group $G'\subset G$. To prove this, we need some
preliminary results which may be of interest in their own right.

First we recall what is a reduction of structure group of a $G$-Higgs bundle. If $(E_G,\varphi)$ is a $G$-Higgs bundle, and $G'$ is a reductive subgroup of $G$, then a \emph{reduction of structure group} of $(E_G,\varphi)$ is a $G'$-Higgs bundle $(E_{G'},\varphi')$ such that $E_{G'}\hookrightarrow E_G$ is a holomorphic reduction of structure group of $E_G$ to the principal $G'$-bundle $E_G'$, and such that $\varphi'$ maps to $\varphi$ under the embedding $\ad(E_{G'})\otimes K\hookrightarrow\ad(E_G)\otimes K$.

\begin{proposition}\label{prop:stable reduction}
Let $(E_G,\varphi)$ be a stable $G$-Higgs bundle which admits a reduction to a $G'$-Higgs bundle for some complex reductive subgroup $G'\subset G$. Let $(E_{G'},\varphi')$ be the $G'$-Higgs bundle obtained by the reduction of the $G$-Higgs bundle $(E_G,\varphi)$. Then $(E_{G'},\varphi')$ is stable as a $G'$-Higgs bundle.
\end{proposition}
\proof Let $H'\subset G'$ be a maximal compact subgroup. Then $\lieg'\subset\lieg$ and $\lieh'\subset\lieh$.
Let $s\in i\lieh'$. Let $P'_s$ be the parabolic subgroup of $G'$ associated to $s$ as defined in \eqref{eq:parabolicPs}. Since also $s\in i\lieh$, it defines a parabolic subgroup $P_s$ of $G$ such that $P'_s\subset P_s$. Now, take a reduction $\sigma'\in H^0(E_{G'}(G'/P'_s))$ and denote by $E_{\sigma'}$ be the corresponding $P'_s$-bundle. Given $\sigma'$ and the reduction of $E_G$ to $E_{G'}$ one naturally obtains an induced reduction $\sigma\in H^0(E_G(G/P_s))$ of $E_G$ to a principal $P_s$-bundle $E_\sigma$, by extending the structure group of $E_{\sigma'}$ through the inclusion $P'_s\hookrightarrow P_s$. 
Moreover, if $\sigma'$ is such that $\varphi'\in H^0(\ad(E_{\sigma'})\otimes K)$, then $\sigma$ is such that $\varphi\in H^0(\ad(E_\sigma)\otimes K)$.
Also, any antidominant character $\chi':P'_s\to\CC^*$ of $P'_s$ gives naturally rise to an antidominant character $\chi:P_s\to\CC^*$ of $P_s$, just by extending $\chi'$ to $P_s$ by $1$, and clearly $\deg_{\chi'}(E_{\sigma'})=\deg_\chi(E_\sigma)$.
 So, from Definition \ref{def:semipoly} of stability, we conclude that if the $G'$-Higgs bundle $(E_{G'},\varphi')$ is not stable, then $(E_G,\varphi)$ is not stable.
\endproof

\begin{proposition}\label{prop:every element of Aut is semisimple}
Let $(E_G,\varphi)$ be a stable $G$-Higgs bundle. Then every element of the Lie group $\Aut(E_G,\varphi)$ is semisimple.
\end{proposition}
\proof
The main point is given by Proposition 2.14 of \cite{garcia-prada-gothen-mundet:2008} which says that the stability of $(E_G,\varphi)$ implies that every element of the infinitesimal automorphism space $\aut(E_G,\varphi)$ is semisimple. Since $\aut(E_G,\varphi)$ is the Lie algebra of $\Aut(E_G,\varphi)$, it follows that every element of $\Aut(E_G,\varphi)_0$ --- the identity component of $\Aut(E_G,\varphi)$ --- is also semisimple. Now, take the projection morphism onto the group of connected components $p:\Aut(E_G,\varphi)\to\pi_0(\Aut(E_G,\varphi))$ and let any $g\in\Aut(E_G,\varphi)$. Then $g=g_sg_u$, where $g_s$ is its semisimple part and $g_u$ its unipotent part. Since $p$ is a morphism it preserves the semisimple and unipotent parts, and since every element of the group $\pi_0(\Aut(E_G,\varphi))$ is semisimple (because it is finite), then $p(g_u)=0$. Hence $g_u\in\Aut(E_G,\varphi)_0$. But we already know that every element in $\Aut(E_G,\varphi)_0$ is semisimple, so $g_u=0$, hence $g=g_s$ is semisimple.
\endproof

From this result we obtain the following corollary, where again we are identifying the elements of $\Aut(E_G,\varphi)$ with elements of the $G$, up to conjugation.
Given one such element $g\in\Aut(E_G,\varphi)$, let $Z_G(g)$ denote the
centralizer in $G$ of $g$. If we choose another representative $hgh^{-1}$ in
$G$ ($h\in G$) for the automorphism defined by $g$ then, since
$Z_G(hgh^{-1})\cong hZ_G(g)h^{-1}$, we see that $Z_G(g)$ is defined up to
conjugation. We have the following.

\begin{corollary}\label{cor:centralizer reductive}
Let $(E_G,\varphi)$ be a stable $G$-Higgs bundle and let 
$g\in\Aut(E_G,\varphi)$. Then $Z_G(g)$ is reductive.
\end{corollary}
\proof
Proving that the complex Lie group $Z_G(g)$ is reductive is (by definition) equivalent to proving that it is the complexification of a compact Lie group. This centralizer is the same as the subgroup of $G$ of fixed points of the inner automorphism of $G$ given by
$\Int(g)(h)=ghg^{-1}$: $$Z_G(g)=\{h\in G\st \Int(g)(h)=h\}.$$
This inner automorphism defines an automorphism of $G_0$, so $\Int(g)\in\Aut(G_0)$. By Proposition \ref{prop:every element of Aut is semisimple}, $g$ is a semisimple element of $G$, hence $\Int(g)$ is also a semisimple element of $\Aut(G_0)$. Let $S$ be the torus in $\Aut(G_0)$ generated by $\Int(g)$, so that the subgroup of $G_0$ fixed by $S$ (or equivalently by $\Int(g)$) is reductive, according to \cite[Proposition 3.6, page 107]{gorbatsevich-onishchik-vinberg:1993}. Denote this reductive group by $Z_{G_0}(g)$.

Now, $Z_G(g)$ is a finite extension of $Z_{G_0}(g)$: just consider the short exact sequence
\begin{equation}\label{eq:finite ext}
0\to Z_{G_0}(g)\to Z_G(g)\to\pi_0(G)\to 0,
\end{equation} where we are taking the restriction to $Z_G(g)$ of the projection $G\to\pi_0(G)$. Let $H$ be a maximal compact subgroup of $Z_G(g)$ (see Theorem 14.1.3 of \cite{hilgert-neeb:2012}). Then $H$ intersects all the components of $Z_G(g)$. So $H\cap Z_{G_0}(g)$ is a maximal compact subgroup of $Z_{G_0}(g)$, whose complexification is precisely $Z_{G_0}(g)$ because we know that $Z_{G_0}(g)$ is reductive. In particular the Lie algebra of $Z_{G_0}(g)$ is the complexification of the Lie algebra of $H\cap Z_{G_0}(g)$ which, by \eqref{eq:finite ext}, is equivalent to say that the Lie algebra of $Z_G(g)$ is the complexification of the Lie algebra of $H$. We can thus apply Proposition 15.2.4 of \cite{hilgert-neeb:2012} to conclude that $Z_G(g)=H^\CC$.
\endproof
 
We can now have a description of stable $G$-Higgs bundles which are not simple. Some of the arguments used in the proof of the following theorem are based on similar ones used in \cite{biswas-parameswaran:2006}.

\begin{theorem}\label{thm:stable and not simple principal bundles}
Let $(E_G,\varphi)$ be a stable $G$-Higgs bundle which is not simple. Then there is a complex reductive subgroup $G'\subset G$ and a reduction of $(E_G,\varphi)$ to a $G'$-Higgs bundle, which is stable and simple.
\end{theorem}
\proof
Since $(E_G,\varphi)$ is not simple, there is some $f\in\Aut(E_G,\varphi)$ which does not belong to $Z(G)$.
Since it is stable, then \cite[Proposition 3.14]{garcia-prada-gothen-mundet:2008} $\aut(E_G,\varphi)\cong\liez$, where $\aut(E_G,\varphi)$ is defined in \eqref{eq:infinit autom}. Notice that if $\varphi=0$, then $E_G$ is stable (cf. Remark \ref{rm:(semi)stability of bundles}) and also $\aut(E_G)\cong\liez$, by Proposition 3.2 of  \cite{ramanathan:1975}. Thus, in any case,
$$\Aut(E_G,\varphi)_0\cong Z(G)_0,$$ $Z(G)_0$ denoting the connected component of $Z(G)$ containing the identity. Moreover, since $\Aut(E_G,\varphi)$ is an algebraic group, it has finitely many connected components, hence the quotient
\begin{equation}\label{eq:Q=Aut/Z}
Q=\Aut(E_G,\varphi)/Z(G)
\end{equation}
is a finite group (note that $Z(G)$ is a normal subgroup of $\Aut(E_G,\varphi)$).
Write
\begin{equation}\label{eq:Q finite}
Q=\{[f_1],\dots,[f_k]\}
\end{equation}
where $f_i\in\Aut(E_G,\varphi)$.

Consider the class $[f_1]$. From this class we shall obtain a reduction of structure group of $E$.
Recall that $f_1$ is a section of $\Ad(E_G)=E_G\times_GG$ and that, by Remark \ref{rem:identif-Aut-G}, the automorphism $f_1$ determines a closure of the orbit of an element $g_1\in G$, under conjugation. However, by Proposition \ref{prop:every element of Aut is semisimple}, $f_1$ is semisimple, hence so is $g_1$. Therefore its orbit by conjugation is closed, so $f_1$ is identified with a conjugacy class $[g_1]$ of $G$. 
Now, let $q:E_G\times G\to\Ad(E_G)$ be the quotient map, and define $$\mathcal I=q^{-1}(f_1(X))\subset E_G\times G.$$
If $p_G:E_G\times G\to G$ and $p_X:E_G\times G\to X$ are the natural projections then, as $f_1(X)=[g_1]$, we have that $$(p_X\times p_G)(\mathcal I)=X\times(G\cdot g_1),$$ where $G\cdot g_1$ denotes the orbit of $g_1$ under the action of $G$ on itself by conjugation. Now, let $\widehat{\mathcal I}\subset\mathcal I$ be given as 
$$\widehat{\mathcal I}=(p_X\times p_G)|_{\mathcal I}^{-1}(X\times\{g_1\}).$$
Let $p_{E_G}:E_G\times G\to E_G$ be the projection, and define
\begin{equation}\label{def:reduction1}
E_1=p_{E_G}|_{\mathcal I}(\widehat{\mathcal I})\subset E_G.
\end{equation}
The restriction to $E_1$ of the projection $\pi:E_G\to X$ is holomorphic and surjective, and the centralizer $Z_{G}(g_1)$ of $g_1$ in $G$ acts transitively on the fibres of $\pi|_{E_1}$. Hence, $E_1$ is a subbundle of $E_G$ whose structure group is $Z_{G}(g_1)$, that is, $E_1$ is a reduction of structure group of $E_G$ to $Z_G(g_1)\subset G$.

Recall now that we have made two choices. Let us see what is the dependence of our reduction on these choices. First, we have chosen a representative $f_1\in\Aut(E_G,\varphi)$ of the class $[f_1]\in Q$, where $Q$ is given by (\ref{eq:Q=Aut/Z}). If $f_1z$ is another representative, with $z\in Z(G)$, then, as above, $f_1z$ will give rise to the conjugacy class $[g_1z]$ of $G$. However, doing the same construction as in the previous paragraph with $g_1$ replaced by $g_1z$, and noticing that $Z_{G}(g_1z)=Z_{G}(g_1)$, one obtains the same reduction of structure group of $E_G$ to the $Z_{G}(g_1)$-bundle $E_1$ as before.
On the other hand, given the class $f_1$, we can choose a different representative $gg_1g^{-1}$ ($g\in G$) of the class $[g_1]$. This will give rise to a reduction of structure group $E'_1$ of $E_G$ to $Z_{G}(gg_1g^{-1})\cong gZ_{G}(g_1)g^{-1}$.
The reduction of structure group of $E_G$ to $Z_{G}(g_1)$ is therefore well-defined, up to conjugation.

Let $\liez_\lieg(g_1)$ be the Lie algebra of $Z_{G}(g_1)$. Then we have $$\liez_\lieg(g_1)=\{v\in \lieg\st\Ad(g_1)(v)=v\}$$ and the adjoint representation restricts to $\Ad_1:Z_{G}(g_1)\to\GL(\liez_\lieg(g_1))$. Moreover, $\varphi\in H^0(\ad(E_1)\otimes K)$, because $g_1$ is an automorphism of $(E_G,\varphi)$ (cf.  \eqref{eq:autom}). Write $\varphi_1$ for $\varphi$ in $H^0(\ad(E_1)\otimes K)$.
Then $(E_1,\varphi_1)$ is a $Z_{G}(g_1)$-Higgs bundle.

Now we iterate this procedure. 
Since $Q$ is finite, this process will end and we obtain a reduction of $(E_G,\varphi)$ to a $G'$-Higgs bundle $(E_k,\varphi_k)$, where $G'=Z_{G}(g_1,\ldots,g_k)$.

The fact that $G'$ is a complex reductive group follows from Corollary \ref{cor:centralizer reductive}. 
The stability of $(E_k,\varphi_k)$ as a holomorphic $G'$-Higgs bundle follows from the stability of $(E_G,\varphi)$ as a $G$-Higgs bundle and from Proposition \ref{prop:stable reduction}.
Finally, since $g_i\in Z(G')$, for every $i=1,\ldots,k$, we conclude that $\Aut(E_k,\varphi_k)=Z(G')$, thus $(E_k,\varphi_k)$ is simple as a $G'$-Higgs bundle.
\endproof

\begin{remark}
As far as we know, Theorem \ref{thm:stable and not simple principal bundles} 
is not in the literature even for principal bundles --- our proof also holds 
in that case by considering $\varphi=0$. 
\end{remark}

\subsection{Strictly polystable Higgs bundles}

Also strictly polystable $G$-Higgs bundles correspond to singularities of
$\M_G$, which this time may be more ``serious'' than those of orbifold
type. However, as we now recall, for such $G$-Higgs bundles there is also a
reduction of structure group such Higgs bundle is stable for the new
group. Indeed, the following result is a particular case for the existence of
a Jordan-H\"older reduction of semistable $G$-Higgs bundles, which is unique 
up to isomorphism (look at \cite{garcia-prada-gothen-mundet:2008,otero:2010}
for two different proofs).

\begin{proposition}\label{prop:strictly polystable}
A polystable $G$-Higgs bundle admits a reduction to a stable $G'$-Higgs bundle, where $G'\subset G$ is a complex reductive subgroup.
\end{proposition}
\proof
If the $G$-Higgs bundle is stable, there is nothing to prove. Assume hence that it is strictly polystable. There is then some parabolic subgroup $P\subset G$, some antidominant character $\chi$ of $P$ and some reduction of structure group $\sigma$ of $E_G$ to $P$ such that $\varphi\in H^0(\ad(E_\sigma)\otimes K)$ such that $\deg_{\chi}(E_\sigma)=0$. Moreover, there is a further holomorphic reduction of structure group $\sigma_L$ of $E_\sigma$ to an principal $L$-bundle where $L$ is the Levi subgroup of $P$ such that $\varphi\in H^0(\ad(E_{\sigma_L})\otimes K)$.
The $L$-Higgs bundle $(E_{\sigma_L},\varphi)$ is polystable (this is proved in the same way as the proof of Proposition \ref{prop:stable reduction}). If it is stable, we are done. If not, we iterate this procedure, which will eventually end, yielding a $G'$-Higgs bundle which must be stable.
\endproof

 \section{The Hitchin function and the subvarieties of local minima}\label{section:Hitchin function}

\subsection{The Hitchin proper function}
In this section we develop all the machinery needed to count the number of connected components of the moduli space of $G$-Higgs bundles $\M_G$, defined in \eqref{eq:modulispace}.

A first division of $\M_G$ into connected is given by the topological class of the Higgs bundles $(E_G,\varphi)$.
So, the number of connected components of $\M_G$ is bounded below by the number of topological classes of $G$-Higgs bundles on the surface $X$, which can be obtain by Theorem \ref{thm:topclass}. If $G$ is connected the lower bound is the cardinal of $\pi_1G$, which may be infinite. So, the real interest is to determine the connected components of $\M_G(c)$, for each class $c$. 

Fix one such class $c$.
In order to study $\pi_0(\M_G(c))$, we use the method introduced by Hitchin in \cite{hitchin:1987}, which uses the non-negative real valued function given by the $L^2$-norm of the Higgs field. When $\M_G(c)$ is smooth, this function is a perfect Morse-Bott function, thus it is clearly a useful tool for the study of the topology of $\M_G(c)$. But even when $\M_G(c)$ is not smooth (which is the large majority of the cases), this function is still proper, so through the study of the connected components of the subvarieties of its local minima, one can draw conclusions about the connected components of $\M_G(c)$. This is by now a standard method, which has been used systematically to study the connected components, and other topological information, of $\M_G(c)$ for many classes of $G$ (see for example \cite{bradlow-garcia-prada-gothen:2005} and the references therein). Our aim is to perform this study from an intrinsic point of view, i.e., without specifying the group $G$.

Consider the real function $f:\M_G(c)\to\RR$ defined as
\begin{equation}\label{eq:Hitfunc}
f(E_G,\varphi)=\|\varphi\|_{L^2}^2=\int_{X}|\varphi|^2\mathrm{dvol}=\int_{X}B(\varphi,\tau_h(\varphi))\mathrm{dvol}.
\end{equation}
where $B$ is a non-degenerate quadratic form on $\lieg$, extending the Killing form on $\lieg_{ss}=[\lieg,\lieg]$ and $\tau_h$ is defined before Theorem \ref{thm:Hit-Kob} and which depends on the metric $h$  which provides the Hitchin-Kobayashi correspondence, that is, the one which is a solution to Hitchin equations --- see Theorem \ref{thm:Hit-Kob}.
This function $f$ is known as the \emph{Hitchin function}.
Using the Uhlenbeck weak compactness theorem, one can prove \cite{hitchin:1987} that $f$ is proper and therefore attains a minimum on each closed subspace $\M'(c)$ of $\M_G(c)$. The next result relates the connectedness of the subspaces of $\M'(c)$ of local minima with the connectedness of $\M'(c)$ itself.

\begin{proposition}\label{proper}
Let $\M'(c)\subseteq\M_G(c)$ be a closed subspace and let $\N'\subset\M'(c)$ be the subspace of local minima of $f$ on $\M'(c)$. If $\N'$ is connected then so is $\M'(c)$.
\end{proposition}

The idea is then to have a detailed description of the subspace of local minima of $f$, enough to draw conclusions about its connectedness. Since this method has been already applied for several cases (see for instance \cite{hitchin:1992,bradlow-garcia-prada-gothen:2013}), we will only sketch it.

The strategy for studying the connectedness of $\M_G(c)$ implies that we resort on a separated approach for the following three disjoint locus of $\M_G(c)$:
\begin{enumerate}
\item[(1)] locus of stable and simple $G$-Higgs bundles;
\item[(2)] locus of stable but not simple $G$-Higgs bundles;
\item[(3)] locus of strictly polystable $G$-Higgs bundles. 
\end{enumerate}

Proposition \ref{prop:stable and simple imply smooth} says that (1) is included in the smooth locus of $\M_G(c)$.

\subsection{Description of the local minima}

 In this section we only consider stable and simple $G$-Higgs bundles.
 A very useful feature of the moduli space $\M_G(c)$ is that it carries a $\CC^*$-action 
\begin{equation}\label{eq:Caction}
\lambda\cdot(E_G,\varphi)=(E_G,\lambda\varphi).
\end{equation}
By considering the moduli space of gauge equivalence classes of solutions to
the Hitchin equations and the restriction of the $\CC^*$-action to an
$S^1$-action, one concludes that a point of $\M_G(c)$ represented by a stable
and simple $G$-Higgs bundle is a critical point of $f$ if and only if it is a
fixed point of the $\CC^*$-action; this is a consequence of the fact that $f$ is a moment map for the $S^1$-action (cf. \cite{hitchin:1987}). Higgs bundles with vanishing Higgs field are obvious fixed points (and global minima of $f$) and the following result (see \cite{hitchin:1987,simpson:1992}) provides a description of the other fixed points and, consequently, of the critical points of $f$.

\begin{proposition}\label{prop:fixed point locus}
A stable and simple $G$-Higgs bundle $(E_G,\varphi)$ is fixed under the $\CC^*$-action \eqref{eq:Caction} if and only if there exists a semisimple element $\psi\in H^0(E_G\times_H\lieh)$ such that there is a decomposition $$\ad(E_G)=\bigoplus_{k=-k_M}^{k_M}\ad(E_G)_k$$ into eigenbundles of $\ad(E_G)$ under the adjoint action $\ad(\psi):\ad(E_G)\to\ad(E_G)$. Here, $k_M$ is positive integer and, for each $k$, $\ad(E_G)_k$ is the $ ik$-eigenbundle. Furthermore, $\ad(\psi)(\varphi)= i\varphi$, so that $\varphi\in H^0(\ad(E_G)_1\otimes K).$
\end{proposition}

If $(E_G,\varphi)$ is a fixed point of the $\CC^*$-action, we can then consider an induced decomposition of the complex $C^\bullet(E_G,\varphi)$, defined in \eqref{eq:deformation complex}, as follows: $$C^\bullet(E_G,\varphi)=\bigoplus_{k=-k_M}^{k_M}C^\bullet(E_G,\varphi)_k,$$
where $C^\bullet(E_G,\varphi)_k$ is the subcomplex of $C^\bullet(E_G,\varphi)$ defined by
\begin{equation}\label{eq:subcomplex}
C^\bullet(E_G,\varphi)_k:\ad(E_G)_k\xrightarrow{\ad(\varphi)_k}\ad(E_G)_{k+1}\otimes K
\end{equation} where we define $\ad(\varphi)_k=\ad(\varphi)|_{\ad(E_G)_k}$. We say that $C^\bullet(E_G,\varphi)_k$ is the subcomplex of \emph{weight} $k$.
In turn, this yields a decomposition $$\mathbb{H}^1(C^\bullet(E_G,\varphi))=\bigoplus_{k=-k_M}^{k_M}\mathbb{H}^1(C^\bullet(E_G,\varphi)_k)$$ of the tangent space of $\M_G(c)$ at $(E_G,\varphi)$.

The following result is fundamental for the description of the smooth local minima of $f$, among the critical points which have just been described. It is a consequence of the fact that, for each $k$, the subspace $\mathbb{H}^1(C^\bullet(E_G,\varphi)_k)$ is the $(-k)$-eigenspace of the hessian of $f$ at $(E_G,\varphi)$. This is basically  \cite[Lemma 3.11]{bradlow-garcia-prada-gothen:2008}. Although the proof in that paper is for $\GL(n,\CC)$-Higgs bundles, the same argument works in the general setting of $G$-Higgs bundles: the key facts are that for a stable $G$-Higgs bundle, $(E_G,\varphi)$, the Higgs vector bundle $(\ad(E_G),\ad(\varphi))$ is semistable, and that there is a natural $\ad$-invariant isomorphism $\ad(E_G)\cong\ad(E_G)^*$ given by the invariant pairing $B$ on $\lieg$, extending the Killing form on $\lieg_{ss}$ --- see the definition of the Hitchin function in \eqref{eq:Hitfunc}.

\begin{proposition}\label{prop:ad isom}
Let $(E_G,\varphi)\in\M_G(c)$ be a stable and simple critical point of $f$. Then
$(E_G,\varphi)$ is a local minimum if and only if either $\varphi=0$ or $\ad(\varphi)_k$ in \eqref{eq:subcomplex} is an isomorphism for all $k\geq 1$.
\end{proposition}

We can now describe the stable and simple local minima of the Hitchin function $f$.

\begin{proposition}\label{prop:smooth minimum}
Let $(E_G,\varphi)$ be a stable and simple $G$-Higgs bundle, which is a critical point of $f$. Then $(E_G,\varphi)$ represents a local minimum if and only if $\varphi=0$.
\end{proposition}
\proof Let $(E_G,\varphi)$ be a local minimum of $f$ and suppose that $\varphi\neq 0$. Consider the complex \eqref{eq:subcomplex}.
Consider the highest possible weight $k_M$, i.e the highest weight such that $\ad(E_G)_{k_M}\neq 0$.
As $\varphi\neq 0$, then $k_M\geq 1$, so by Proposition \ref{prop:ad isom}, $$\ad(E_G)_{k_M}\cong \ad(E_G)_{k_M+1}\otimes K=0,$$ a contradiction.
\endproof

Now we consider $G$-Higgs bundles which are stable but not simple. In principle, these do not correspond to smooth points of the moduli space, so our analysis of the local minima of this kind must be carried out in a different way, because Proposition \ref{prop:ad isom} does not apply.
However, given their description by Theorem \ref{thm:stable and not simple principal bundles}, we now achieve easily the goal of describing the local minima of the Hitchin function \eqref{eq:Hitfunc} on the stable locus of the moduli space $\M_G(c)$.

\begin{proposition}\label{prop:stable minimum}
Let $(E_G,\varphi)$ be a stable $G$-Higgs bundle. Then $(E_G,\varphi)$ represents a local minimum of the Hitchin function if and only if $\varphi=0$.
\end{proposition}
\proof
The stable and simple case is the content of Proposition \ref{prop:smooth minimum}. Assuming hence that $(E_G,\varphi)$ is not simple, Theorem \ref{thm:stable and not simple principal bundles} assures the existence of a complex reductive subgroup $G'\subset G$ and a reduction of $(E_G,\varphi)$ to a $G'$-Higgs bundle $(E_{G'},\varphi')$, which is stable and simple. On the corresponding moduli space $\M_{G'}$, $(E_{G'},\varphi')$ must be a local minima of the restriction of Hitchin function to $\M_{G'}$. Proposition \ref{prop:smooth minimum} implies that $\varphi'=0$, thus $\varphi=0$.
\endproof

 The remaining case to be considered is the locus of strictly polystable $G$-Higgs bundles.
From Proposition \ref{prop:strictly polystable}, we achieve the description of the subvariety of local minima of $f$ in $\M_G(c)$. Let $\N_G$ denote the moduli space $\N_G$ of polystable principal $G$-bundles over $X$, and let $\N_G(c)$ be the subspace of $\N_G$ given by those $G$-bundles with topological class determined by $c$.

\begin{theorem}\label{thm:local minimum}
The subvariety of local minima of the Hitchin function over $\M_G(c)$ is isomorphic to $\N_G(c)$.
\end{theorem}
\proof
Let $(E_G,\varphi)$ be a polystable $G$-Higgs bundle which is a local minimum of the Hitchin function. If it is stable, Proposition \ref{prop:stable minimum} states that $\varphi$ must vanish. If it is strictly polystable then, using Proposition \ref{prop:strictly polystable} and along the same lines as the proof of Proposition \ref{prop:stable minimum}, one concludes that $\varphi=0$ as well.
So the local minima of $f$ is the subvariety of $\M_G(c)$ given by those $(E_G,0)$, which is isomorphic to $\N_G(c)$.
\endproof

\section{Connected components}\label{compM}

In \cite{ramanathan:1996b}, Ramanathan has shown that if $G$ is a connected reductive complex Lie group, then $\N_G(c)$ is connected for topological class $c$ in $\pi_1G$. However his arguments readily adapt to the case where $G$ is non-connected. For the benefit of the reader we provide the details.

\begin{proposition}\label{prop:connectedness of moduli of principal bundles}
For a complex reductive Lie group $G$  and any topological class $c$, the moduli space $\N_G(c)$ is connected.
\end{proposition}
\proof Let $E'_G$ and $E''_G$ represent two classes in
$\N_G(c)$. Let $P$ be the underlying $C^\infty$ principal bundle, and let
$\overline{\partial}_{A'}$ and $\overline{\partial}_{A''}$ be the
operators on $P$ defining, respectively, $E'_G$ and $E''_G$ and given by
$H$-connections $A'$ and $A''$, respectively, where $H$ is a maximal compact subgroup of $G$.

Let $\mathbb D$ be an open disc in $\CC$ containing 0 and 1. Consider the
$C^\infty$ principal $G$-bundle
$\mathbb E_G\to\mathbb D\times X$, where $\mathbb E_G=\mathbb D\times P$. Define the
connection form on $\mathbb E_G$ by
$$A_z(v,w)=zA''(w)+(1-z)A'(w)\in\Omega^1(\mathbb E_G,\lieg)$$
where $v$ is tangent to $\mathbb D$ at $z$ and $w$ is tangent to $P$ at some
point $p$. If we consider the holomorphic bundle $E_z$ given by
$\mathbb E_G|_{\{z\}\times X}$ with the holomorphic structure given
by $A_z$, then we have that $E_0\cong E'_G$ and $E_1\cong E''_G$.

As semistability is an open condition with respect to the Zariski
topology, $\mathbb D\setminus D'$ is connected where $D'=\{z\in \mathbb D:E_z\text{ is not
semistable}\}$. Hence $\{E_z\}_{z\in \mathbb D\setminus D'}$ is a connected family of
semistable principal $G$-bundles joining $E_0$ and $E_1$. Since $E_0\cong E'_G$ and
$E_1\cong E''_G$, using the universal property of the coarse moduli space 
$\N_G(c)$ of $G$-principal bundles, we conclude that there is
a connected family in $\N_G(c)$ joining $E'_G$ and
$E''_G$.
\endproof

We can now state our main result.

\begin{theorem}\label{main-theorem}
Let $G$ be a complex reductive Lie group and let $c$ be a topological class of $G$-Higgs bundles. Then $\M_G(c)$ is non-empty and connected. Thus there is a bijection between $\pi_0(\M_G)$ and the set of topological classes of $G$-Higgs bundles which, when $\pi_0G$ is abelian, are given by Theorem \ref{thm:topclass}. In particular, if $G$ is connected, the number of connected components of $\M_G$ equals the cardinal of $\pi_1G$.
\end{theorem}
\proof
Immediate from Corollary \ref{cor:Mc(G)non-empty}, Proposition \ref{proper}, Theorem \ref{thm:local minimum} and Proposition \ref{prop:connectedness of moduli of principal bundles}. 
\endproof

Let $\Gamma$ be the universal central extension of $\pi_1X$. It is the finitely generated group defined as follows:
$$\Gamma=\bigg<a_1,b_1,\ldots,a_g,b_g,\delta\st\prod_{i=1}^g[a_i,b_i]=\delta\text{ and }\delta\text{ central}\bigg>.$$ Clearly one has the extension $0\to\Z\to\Gamma\to\pi_1X\to 0$.
Now, let $$\Gamma_\RR=\RR\times_\Z\Gamma.$$ Then we have an extension $$0\to\RR\to\Gamma_\RR\to\pi_1X\to 0.$$ A representation $\rho:\Gamma_\RR\to G$ is a continuous homomorphism of groups. It is called \emph{central} if $\rho(\RR)\subset Z(G_0)$. In fact, since $\rho$ is continuous and $\RR$ contains the unit element, we must have $\rho(\RR)\subset Z(G_0)_0$. Notice that if $G$ is semisimple then a central representation $\rho:\Gamma_\RR\to G$ is really a representation of $\pi_1X$ in $G$.
For a central representation $\rho$, we must then have $\rho(\delta)\in Z(G_0)_0$, but since $\prod_{i=1}^g[a_i,b_i]=\delta$, we must in fact have $\rho(\delta)\in Z(G_0)_0\cap G_{ss}$ where $G_{ss}$ denotes the semisimple part of $G$. 

Let $\Hom_{\text{cent}}^{\text{red}}(\Gamma_\RR,G)$ the space of reductive representations (i.e. the ones which become completely reducible when composed with $\Ad:G\to\GL(\lieg)$) which are central. The group $G$ acts by conjugation on this space, as $(g\cdot\rho)(\gamma)=g\rho (\gamma)g^{-1}$. The \emph{$G$-character variety of $X$} is the quotient space $$\cR_G=\Hom_{\text{cent}}^{\text{red}}(\Gamma_\RR,G)/G.$$ 

Non-abelian Hodge theory on $X$ establishes a homeomorphism between $\M_G$ and
$\cR_G$
(cf. \cite{hitchin:1987,simpson:1988,simpson:1992,garcia-prada-gothen-mundet:2008}). Let
$\cR_G(c)$ be the subspace of $\cR_G$ whose corresponding projectively flat
$G$-bundle belongs to the topological class $c$. A direct 
consequence of Theorem \ref{main-theorem}  is hence the following.
\begin{theorem}\label{thm:main}
Let $G$ be a complex reductive Lie group and let $c$ be a topological class of $G$-Higgs bundles. Then $\cR_G(c)$ is non-empty and connected. Thus there is a bijection between $\pi_0(\cR_G)$ and the set of topological classes of $G$-Higgs bundles which, when $\pi_0G$ is abelian, are given by Theorem \ref{thm:topclass}. In particular, if $G$ is connected, the number of connected components of $\cR_G$ equals the cardinal of $\pi_1G$.
\end{theorem}

\begin{remark}\label{rmk:knowncases}\mbox{}
\begin{enumerate}
\item
In the case of $G$ complex, reductive and connected, this result also follows from the work \cite{donagi-pantev:2012} of Donagi and Pantev. The methods used there are however completely different than the ones used in this paper. There the Hitchin system is the main tool and the authors prove that the generic fibre of the Hitchin map has precisely $\pi_1G$ connected components. From this the conclusion that $\pi_0(\M_G)=\pi_1G$ follows. On the other hand, it seems to us that our method is more suitable for considering a generalisation of this intrinsic study of $\pi_0(\M_G)$ for any real reductive Lie group. 

\item If $G$ is connected and semisimple, the main result of this paper has also already been proved by J. Li, in \cite{li:1993}, using different methods, involving the study of flat bundles on Riemann surfaces. Our result provides thus also an alternative proof to Li's theorem in that case. Finally, very recently, Li's result has been generalised, by N-K. Ho and C-C. M. Liu, to include representations of the fundamental group of $X$ into a complex, connected and reductive Lie group --- see the Appendix of \cite{lawton-ramras:2014}. Notice that they study the components of $\Hom(\pi_1X,G)$, but these are in bijection to the ones of $\Hom^\mathrm{red}(\pi_1X,G)$, since, by
Theorem 8 of \cite{garcia-prada-mundet:2004}, every representation from $\pi_1X$ to $G$ can be deformed to a reductive one. Their method is a direct approach, using representations and using also Li's theorem. However, the representations considered there do not include all the possible topological types (for that, one needs to consider representations of $\Gamma_\RR$), so our result is a generalisation of theirs, even in the case when $G$ is connected.
\end{enumerate}
\end{remark}

\noindent 
\textbf{O. Garc{\'\i}a-Prada}\\
Instituto de Ciencias Matem\'aticas,  CSIC \\ 
Nicol\'as Cabrera, 13--15, 28049 Madrid, Spain\\
email: oscar.garcia-prada@icmat.es

\vspace{1cm}

\noindent
      \textbf{André Oliveira} \\
      Centro de Matemática da Universidade do Porto, CMUP\\
      Faculdade de Ciências, Universidade do Porto\\
      Rua do Campo Alegre 687, 4169-007 Porto, Portugal\\ 
      \url{www.fc.up.pt}\\
      email: andre.oliveira@fc.up.pt

\vspace{.2cm}
\noindent
\textit{On leave from:}\\
 Departamento de Matemática, Universidade de Trás-os-Montes e Alto Douro, UTAD \\
Quinta dos Prados, 5000-911 Vila Real, Portugal\\ 
\url{www.utad.pt}\\
email: agoliv@utad.pt
\end{document}